\documentclass[10pt]{amsart}
\usepackage{pslatex}

\usepackage{upref,amsxtra,amssymb,amscd}
\usepackage{varioref}
\usepackage{verbatim}
\usepackage{epsfig}

\usepackage{epic,eepic,eucal}
\usepackage{enumerate}

\newcommand{\NN}{{\mathbb N}}

\newcommand{\TT}{{\mathbb T}}
\newcommand{\jk}{{\mathbb S}}
\newcommand{\ZZ}{{\mathbb Z}}

\newcommand{\DD}{{\mathbb D}}

\def\carre{ \hfill $\Box$    }

\newtheorem{theo}{Theorem}[section]
\newtheorem{prop}[theo]{Proposition}

\newtheorem{coro}[theo]{Corollary}

\theoremstyle{definition}

\theoremstyle{remark}

\newtheorem{rema}[theo]{Remark}

\numberwithin{equation}{section}

\begin{document}
\title{Non uniform hyperbolicity and elliptic dynamics}
\author[Bassam Fayad]{Bassam Fayad}
\address{Bassam  Fayad, LAGA, Universit\'e Paris 13, Villetaneuse, }
\email{fayadb@math.univ-paris13.fr}


\begin{abstract}

We present some constructions that are merely the fruit of combining recent results from two areas of smooth dynamics: nonuniformly hyperbolic systems and elliptic constructions.

Our main construction gives a positive answer to a problem raised in \cite{burnsdolpesin}: There exists a volume preserving diffeomorphism which displays a nonuniformly hyperbolic ergodic component that is open and dense but has arbitrarily small measure.  The construction is based on the following techniques:

\begin{itemize}

\item Special versions of the  {\sl approximation by conjugation} method introduced in \cite{AK} such as the one used in \cite{FK} to produce   transitive diffeomorphisms with a zero measure set of points having a dense orbit.

\item Perturbation techniques applied to partially hyperbolic systems that lead to nonuniformly hyperbolic and accessible systems \cite{shubwil,dolpesin,dolwil}.


\end{itemize} 

\end{abstract}

\maketitle

\section{Introduction}

\subsection{} \label{111}

The ergodicity of the geodesic flow on a compact rank 1 manifold is a major open problem.  In the special case of a compact negatively curved space, Anosov proved the ergodicity of the geodesic flow using its uniform hyperbolic structure. In the general case  it is known that the flow is ergodic and nonuniformly hyperbolic on an open and dense set of the phase space but it is still unknown whether this open and dense ergodic component has to be of full measure or not.

Perhaps inspired by this problem, the authors of \cite{dolhupesin} pose the following question in the general context of smooth dynamics: {\sl Is there a volume preserving diffeomorphism which has nonzero Lyapunov exponents on an open {\rm (mod 0)} and dense set $U$ which has positive but not full measure? Is there a volume preserving diffeomorphism with the above property such that $f|U$ is ergodic? }

We give a positive  answer to this problem with a construction based on the following scheme: Start  with a map on the torus $\TT^4$,  $F= A \times {\rm Id}_{\TT^2}$, where $A$ is a hyperbolic automorphism of the torus $\TT^2$. It was discovered in \cite{shubwil} that it is possible to perturb $F$ so that the sum of the integrated  Lyapunov exponents in the center becomes negative. Given a small square  $ C \subset \TT^2$ it is possible, by a special version of the Shub and Wilkinson perturbation accomplished in \cite{dolpesin}, to find a perturbation of $F$ with support on $\TT^2 \times C$ that renders the integrated central exponents of the perturbed map both strictly negative (Cf. \S \ref{neg}). This is an open condition, so applying the result by Dolgopyat and Wilkinson asserting the $C^1$ density of accessible systems (Cf. \S \ref{accessibility}) among partially hyperbolic ones we can further perturb  on $\TT^2 \times C$ to get stable-accessibility on $\TT^2 \times C'$ where $C'$ is a square strictly included in $C$. Because the systems we deal with are small $C^1$ perturbations of $F$ the Pugh and Shub conditions of central bunching and dynamical coherence are satisfied so that an accessibility class is included in a single ergodic component, hence $\TT^2 \times C'$ is included in an open ergodic component where almost every point has strictly negative    central exponents, and this property will persist under small perturbation. To extend the latter ergodic component outside $\TT^2 \times C$ we replace $F$ in the beginning of the above procedure by $A \times T$ where $T$ is a small perturbation of ${\rm Id}_{\TT^2}$  that is transitive but has a large measure set $K$ of points that do not visit $C$. As a result we get a diffeomorphism with an open and dense ergodic component where the central exponents are strictly negative which corresponds to the extended accessibility class of a point in $\TT^2 \times C'$, and an invariant nowhere dense closed set $\TT^2 \times K$ of arbitrarily large measure (depending on how small we take $C$) where the central exponents are equal to zero.

\subsection{} A modification in the last step of the above scheme is the following: use for  $T$ a volume preserving  ergodic diffeomorphism of the disc $\DD^2$ obtained {\sl via } successive conjugations as in \cite{AK}. We get then an  ergodic diffeomorphism of $\TT^2 \times \DD^2$ with an extended accessibility class of full measure and with  strictly negative central exponents. However, an arbitrarily small perturbation of this system, simply obtained by an infinitesimal rotation  applied to $T$ in $\DD^2$, breaks down by Herman's theorem the ergodicity of $T$ near the boundary of the disc and subsequently the accessibility and negative exponents in the central direction for the resulting map of $\TT^2 \times \DD^2$.

\subsection{} Another modification of the construction described in \S \ref{111} is the following: Let $C_1$ and $C_2$ be two disjoint squares in $\TT^2$, then apply the first steps of the construction separately in $\TT^2 \times C_1$ and $\TT^2 \times C_2$. Using the method of successive conjugations one can construct a map $T$ on $\TT^2$ such that $C_1$ and $C_2$ are contained in two open ergodic components  $O_1$ and $O_2$ respectively, such that $\overline {O_1 \cup O_2} = \TT^2$ and $O_1$ and $O_2$ are separated by a pseudo-circle. Since this map $T$ can be obtained arbitrarily close to Identity it can be used in the last step of the construction yielding a non uniformly hyperbolic diffeomorphism of $\TT^4$ with two ergodic components $\TT^2 \times O_1$ and $\TT^2 \times O_2$.

\subsection{} The above constructions present intricate combinations between elliptic and hyperbolic behavior and illustrate how the use of sophisticated elliptic phenomena in the central direction of partially hyperbolic systems can lead with the available perturbation techniques to new examples of nonuniformly hyperbolic behavior, especially in high dimensions. Since various elliptic constructions have symplectic versions as well as nonconservative ones the nonuniformly hyperbolic examples obtained using them can be extended to Hamiltonian systems as well as to dissipative ones. In \cite{hertz}, F. Hertz proves stable ergodicity of non-Anosov ergodic linear automorphism using the Diophantine elliptic behavior in the central direction.  This celebrated work pushes forward the Pugh and Shub program on the theme of ``a little hyperbolicity goes a long way in guaranteeing stably ergodic behavior'' \cite{pughshub5}. We advocate that the use of Liouvillean behavior in the central direction can be useful in producing  complementary examples to the theory.

\subsection{} {\sl Partial hyperbolicity, nonuniform hyperbolicity and local ergodicity.} Let $f$ be a diffeomorphism of a manifold $M$ preserving a Riemannian volume $\mu$. The Lyapunov exponent for a vector $v$ at a point $x \in M$ is the limit:
$$\lambda(x,v) = \limsup_{n \rightarrow +\infty} {1 \over n} \parallel d_x f^n v \parallel.$$ 
  For a fixed $x$, this limit can take at most ${\rm dim} (M)$ values called the Lyapunov exponents of $f$ at $x$ and these values are constant on ergodic components of $f$. We say that a set ${\mathcal O}$ such that $\mu( {\mathcal O}) > 0$ is {\sl nonuniformly hyperbolic} for $f$ if for almost every $x \in {\mathcal O}$ all the Lyapunov exponents of $f$ at $x$ are nonzero.


A diffeomorphism $f$ of a compact smooth Riemannian manifold $M$ is called {\sl partially hyperbolic} if for every $x \in M$ the tangent space at $x$ admits an invariant splitting 
$$T_x M = E^s(x) \oplus E^c(x) \oplus E^u(x) $$
with at least two of the subbundles in the sum nontrivial, and there exists constants $\lambda_s < \lambda_c \leq 1 \leq \lambda'_c < \lambda_u ,$ such that for all $x \in M$ and all $v \in T_xM$ 
\begin{eqnarray*}
v \in E^s(x)  &\Rightarrow& \parallel d_x f (v) \parallel \leq  \lambda_s \parallel v \parallel, \\
v \in E^c(x)  &\Rightarrow&  \lambda_c \parallel v \parallel \leq \parallel d_x f (v) \parallel \leq  \lambda'_c \parallel v \parallel, \\
v \in E^u(x)  &\Rightarrow&  \lambda_u \parallel v \parallel \leq \parallel d_x f (v) \parallel.
\end{eqnarray*}

The bundles  $E^s(x) = E^s_f(x)$ and $E^u(x) = E^u_f(x)$ are called {\sl  stable} and {\sl unstable} bundles respectively, and $E^c(x) = E^c_f(x)$  the {\sl center} bundle of $f$.

In \cite{burnsdolpesin} the authors show that if a volume preserving partially hyperbolic diffeomorphism $f$ has strictly negative central exponents on an invariant set $\mathcal{O}$ of positive measure then every ergodic component of $f | {\mathcal O}$ is open mod (0). This is called a {\sl local ergodicity} phenomenon and is based on the absolute continuity of the local stable manifolds and the local strongly unstable ones as well as on the uniform size of the latter.
If in addition $f$ is assumed to be transitive then ${\mathcal O}$ must be dense and $f| {\mathcal O}$ ergodic.
The question of whether ${\mathcal O}$ has full measure or not naturally arises then, and they answer it positively under the additional assumption of essential accessibility, pointing that in fact only $\epsilon$ accessibility is enough where $\epsilon$ depends on $f$.  Indeed, an argument based on the existence of {\sl hyperbolic times} similar to that of \cite{ABV}  allows them to assert that since $f$ is partially hyperbolic and has negative central exponents, any of its ergodic components contains a ball of uniform size.  
Then $\epsilon$-accessibility is used because it is stable under perturbation and because it guarantees that the orbit of almost every point is $\epsilon$-dense.  This is of course not implied in general by transitivity and smooth volume  preserving transitive diffeomorphism was constructed in \cite{FK} that  have a zero measure set of points with dense orbits.  Based on the latter constructions and on the perturbation techniques of partially hyperbolic systems into nonuniformly hyperbolic accessible systems we prove the following

\begin{theo} For any $\varepsilon >0$ there exists a $C^\infty$ partially hyperbolic diffeomorphism $f$ of the four dimensional torus $\TT^4$  such that 
\begin{enumerate}
\item $f$ preserves the Riemannian volume $\mu$ on $\TT^4$;
\item $f$ has an open and dense nonuniformly hyperbolic ergodic component $\mathcal{ O}$, the central exponents for almost every point in ${\mathcal O}$ being strictly negative;
\item $f$ has a closed invariant set $K$ with $\mu(K) \geq 1 - \varepsilon$ such that  almost every point in $K$ has two strictly positive Lyapunov exponents, two strictly negative and two equal to zero. 
\end{enumerate}
\end{theo}

\begin{rema} It is natural to ask whether it is possible to modify the example in such a way that instead of (3) $f$ becomes nonuniformly hyperbolic on all $M$. The only possibility  would be for $f$ to have on the complementary of $O$ nonzero central exponents of different signs.  With the available techniques however,  it is not clear how can one remove zero exponents on a set that is not open and keep it invariant, and it is a challenging problem in general to find  examples of volume preserving diffeomorphisms having  nonuniformly ergodic components that are not open. \end{rema}

\section{Negative integrated central exponents.} \label{neg}

\subsection{\sl One dimensional central direction} Let $A: \TT^2 \rightarrow \TT^2$ be a linear hyperbolic automorphism. 
Consider the map $F_1 = A \times {\rm Id}_{\jk^1}$ of the three dimensional torus $M_1= \TT^2 \times \jk^1$. Following \cite{shubwil}, it is possible to apply a small perturbation to  $F_1$ in order to get a map with a nonzero integrated central exponent. However, such a perturbation generically gives to the perturbed system the accessibility feature and we end up with an ergodic nonuniformly hyperbolic diffeomorphism. To obtain our example we will have to restrict in the beginning the perturbation to subsets of the type   $\TT^2 \times I$ where $I$ is an interval of the circle as in \cite{dolhupesin}. 

\subsection{\sl Two dimensional central direction.} Moreover, as we explained in the introduction our strategy requires the use in the central direction of a transitive volume preserving map with a large  measure set of points staying away from some small open set.  Since such maps do not exist on the circle we have to consider maps having at least two dimensional central direction. So we start with
$F=A \times {\rm Id}_{\TT^2}$ on the torus $\TT^4$. The additional difficulty in this case is that only the sum of the two exponents in the central direction can be guaranteed to become negative by the Shub and Wilkinson type perturbation and an additional perturbation similar to  the one used in \cite{dolpesin} and based on the Bochi Ma\~n\'e method \cite{bochi}  is necessary to make the exponents in the center close to each other and therefore both negative. To make this additional perturbation possible however, one must guarantee that the central direction does not display a dominated splitting after the Shub Wilkinson type perturbation, which can be done by {\sl saving} a fiber above a fixed point of $A$ where the map  obtained from $F$ is kept equal to identity. It is in the neighborhood of this fiber that the subsequent Bochi Ma\~n\'e type perturbation is applied. Because we will need it in the sequel we will assume that $A$ has at least two fixed point and require that at the end of the two step procedure we have just described the map $\overline{F}$ we obtain remains equal to identity on the fiber above one of the fixed points (Condition (2) below).  The following proposition can be extracted from \cite{dolpesin}

\begin{prop} Let $F= A \times {\rm Id}_{\TT^2}$. For $\varepsilon >0$ let 
$$M^\varepsilon := \TT^2 \times {[0, \varepsilon]}^2.$$
 Then there exists a $C^\infty$  diffeomorphism $h: M \rightarrow M$, that can be chosen arbitrarily close to identity in the $C^1$ topology, such that for $\overline{F} = F \circ h$ we have 

\begin{enumerate}

\item $h = {\rm Id}$ on $M \setminus M^\varepsilon$;

\item $h = {\rm Id}$ on an open set ${\mathcal V}$ such that   $\lbrace 0,0 \rbrace \times \TT^2 \subset {\mathcal V}$;

\item  $\displaystyle{ \int_{M} \log \parallel D \overline{F} | E^c_{ \overline{F}} (w) \parallel  dw = \int_{M^\varepsilon} \log \parallel D  \overline{F} | E^c_{ \overline{F}} (w) \parallel  dw < 0.}$

\end{enumerate} 
Moreover, if we set 

$$A^{\varepsilon} := \TT^2 \times {\left[ {\varepsilon \over 4},  { 3 \varepsilon \over 4} \right]}^2 \subset M^\varepsilon,$$ 
it is possible to choose $h$  such that  

\begin{itemize}

\item[(4)]  $\displaystyle{ \int_{M \setminus A^\varepsilon} \left| \log \parallel D  \overline{F} | E^c_{ \overline{F}} (w) \parallel \right| dw < {1 \over 10}  \left| \int_{A^\varepsilon} \log \parallel D  \overline{F} | E^c_{ \overline{F}} (w) \parallel  dw \right|.}$

\end{itemize} 

\label{p2}

\end{prop}

\section{Accessibility inside $M^\varepsilon$.} \label{accessibility}

A powerful tool to study ergodicity for a partially hyperbolic system is accessibility. Two points $p,q \in M$ are called {\sl accessible } if there are points $p= z_0, z_1, \dots, z_l = q$, $z_i \in M$ such that $z_i \in W^u(z_{i-1})$ or $z_i  \in W^s(z_{i-1})$ for $i=1, \dots, l$. Accessibility is an equivalence relation. The diffeomorphism $f$ is said to be {\sl accessible on a set $M' \subset M$} if $M'$ is included inside a single accessibility class.

We now proceed to the second step of the construction which consists of making $A^\varepsilon$ belong to the same accessibility class. For this we rely on the work \cite{dolwil} that proves the $C^1$ density of accessible systems among partially hyperbolic ones and from which  one can extract the following   

\begin{prop} \label{p3}  Let ${F}$ and $h$ and ${\mathcal V}$ be as in Proposition \ref{p2}. There exists a $C^\infty$ volume preserving diffeomorphism $\tilde{h}$ of $M$, that can be chosen arbitrarily close to identity in the $C^1$ topology,  for which there exists  $\eta >0$ such that for any diffeomorphism $G$ of $M$ satisfying ${\parallel G - F \parallel}_{C^1} \leq \eta$ we have  for $g= \tilde{h} \circ G \circ h$

\begin{enumerate}

\item $\tilde{h} = {\rm Id}$ on $M \setminus M^\varepsilon$;

\item $\tilde{h} = {\rm Id}$ on ${\mathcal V}$;

\item $g$ is accessible  on $A^\varepsilon$;

\item $ \displaystyle{ \int_{A^\varepsilon} \log \parallel D g | E^c_{g} (w) \parallel  dw  < 0}$;

\item $ \displaystyle{ \int_{M \setminus A^\varepsilon} \left| \log \parallel D g | E^c_{g} (w) \parallel  \right| dw  < {1 / 2} \left| \int_{A^\varepsilon} \log \parallel D g | E^c_{g} (w) \parallel dw \right|}$; 

\end{enumerate} 

\end{prop}

\noindent {\sl Proof. } From \cite{dolwil} we can find $\tilde{h}$ satisfying (1) and (2) such that for all $G$ sufficiently close to $F$ we  have (3) Since in addition $\tilde{h}$ can be chosen arbitrarily close to identity in  the $C^1$ topology, (4) and (5) will follow by continuity from (3) and (4) of Proposition \ref{p2}.  \carre

\begin{rema} Since $\tilde{h}$ and $h$ can be chosen arbitrarily close to identity in the $C^1$ topology and $\eta$ arbitrarily small, we can assume that the diffeomorphism $g=\tilde{h} \circ G \circ h$ is a sufficiently small perturbation of $A  \times {\rm Id}_{\TT^2}$ and hence satisfies the conditions of Pugh and Shub (central bunching and dynamical coherence) guaranteeing that an accessibility class is included in a single ergodic component \cite{shubpugh}. In the sequel we will assume that $h$ and $\tilde{h}$ are indeed chosen sufficiently close to identity and $\eta$ sufficiently small so that (3) implies that $A^\varepsilon$ belongs to one ergodic component of  $\tilde{h} \circ G \circ h$. Moreover (4) and (5) then imply that the central Lyapunov exponents on the ergodic component of $A^\varepsilon$ are strictly negative.  \end{rema}

\section{Substituting identity by an elliptic map in the central direction.}

We can summarize the previous steps as follows: There exists $h$ and $\tilde{h}$ with support in $M^\varepsilon$ such that for any $G$ sufficiently close to $F= A \times {\rm Id}_{\TT^2}$ we have that $A^\varepsilon$ is included in a nonuniformly hyperbolic ergodic component for $\tilde{h} \circ G \circ h$. In the case of $ \tilde{h} \circ F \circ h$ the ergodic component containing $A^\varepsilon$ is contained in $M^\varepsilon$ so it has a  measure smaller than $\varepsilon^2$.

We now proceed to the third and last step of our construction in which we want to extend the ergodic component containing $A^\varepsilon$ to an open and dense set while maintaining its measure small. We will do so by taking $G$ equal to $A  \times T$ where $T$ is a transitive map of $\TT^2$ having a large measure set of points staying outside of $ {[0,\varepsilon]}^2$. Such a map was constructed in \cite{FK} using a special version of the method of successive conjugations introduced in \cite{AK}.  We have

\begin{prop} \label{p4} For any $\varepsilon > 0$, there exists a $C^\infty$
diffeomorphism $T$ of $\TT^2$, that can be chosen arbitrarily close to identity in the $C^\infty$ topology, such that 

\begin{enumerate}

\item $T$ preserves the Haar measure $\nu$ on $\TT^2$;

\item $T$ is transitive;

\item There exists a closed set $K$ invariant by $T$ such that $\nu (K) \geq 1- \varepsilon$ and $K \subset \TT^2 \setminus  {[0,\varepsilon]}^2$.

\end{enumerate}

\end{prop}

As a corollary of this Proposition and Proposition \ref{p3} we get for $\tilde{h}$ and $h$ as in the previous sections 

\begin{coro} \label{corr} We can choose $T$ in Proposition \ref{p4} sufficiently close to identity so that the map $f= \tilde{h} \circ (A \times A \times T) \circ h$ satisfies 

\begin{enumerate}

\item $f$ preserves the Riemannian volume on $M= \TT^4$;

\item $A^\varepsilon$ is included in an open and dense ergodic component of $f$ and the Lyapunov exponents in the central direction of $f$ are strictly negative for almost every point in this component; 

\item $\TT^2 \times K$ is invariant by $f$ and $f_{|{\TT^2 \times K}} = {(A \times T)}_{|{\TT^2 \times K}}$, hence $f$ has two strictly positive Lyapunov exponents, two strictly negative and two equal to zero  for almost every point in $\TT^2 \times K$. 

\end{enumerate}

\end{coro}

\noindent {\sl Proof.} The extended accessibility class of a set $B$ is the set of points accessible from a point in $\cup_{n \in \ZZ} f^n (B)$. In our context, the extended accessibility class of a set $B$ is included ${\rm mod} (0)$  inside a single ergodic component. Hence we finish if we prove that the extended accessibility class of $A^\varepsilon$, that we denote by ${\mathcal A}$, is dense. For this we  use the fact that $H = \lbrace 0,0 \rbrace \times \TT^2 $ is invariant by $f$ and that $f_{|H} = T$, hence is transitive. The latter implies that ${\mathcal A} \cap H $ is dense in $H$.

On the other hand, for any $x \in M$ there is an $x' \in H$ such that $x$ is in the accessibility class of $x'$.  For any $\epsilon >0$,  from the density of  ${\mathcal A} \cap H$ in $H$ we can find $\tilde{x}' \in {\mathcal A} \cap H$ and $\tilde{x}$ in the accessibility class of $\tilde{x}'$ (hence $\tilde{x} \in {\mathcal A}$)   such that $d(\tilde{x},x) \leq \epsilon$. The proof of (2) is over. The zero exponents in (3) follow from the fact that the map $T$ has zero exponents at every point of $\TT^2$ because there exists  a sequence $q_n \rightarrow \infty$ such that $T^{q_n}$ converges uniformly to the Identity in the $C^r$ norm for any $r \in \NN$.     \carre

\vspace{0.2cm} 

\noindent {\bf Acknowledgment.} I thank Dmitri Dolgopyat for many useful comments and for helping  me fill a gap in the original construction. I also thank Christian Bonatti, Federico Hertz, Fran\c{c}ois Ledrappier, Yakov Pesin and Marcelo Viana for stimulating conversations. 

\frenchspacing
\bibliographystyle{plain}

\end{document}